\documentclass[12pt,letterpaper]{article}
\usepackage{ifthen,latexsym,amssymb,amsmath}

\renewcommand{\baselinestretch}{1.17}
\newcommand{\C}[1]{{\protect\cal #1}}

\newcommand{\I}[1]{{\mathbb #1}}

\newcommand{\ceil}[1]{\lceil #1\rceil}
\newcommand{\e}{\varepsilon}
\newcommand{\floor}[1]{\lfloor #1\rfloor}

\renewcommand{\mid}{:}

\newcommand{\llabel}[1]{\label{#1}}

\newcommand{\beq}[1]{\begin{equation}\llabel{eq:#1}}
\newcommand{\eeq}{\end{equation}}
\newcommand{\req}[1]{\textrm{(\ref{eq:#1})}}

\newtheorem{theorem}{Theorem}
\newcommand{\bth}[2][nothing]{\ifthenelse{\equal{#1}{nothing}}
 {\begin{theorem}} {\begin{theorem}[#1]}\llabel{th:#2}}

\newtheorem{example}[theorem]{Example}

\newtheorem{lemma}[theorem]{Lemma}
\newcommand{\blm}[2][nothing]{\ifthenelse{\equal{#1}{nothing}}
 {\begin{lemma}} {\begin{lemma}[#1]}\llabel{lm:#2}}

\newtheorem{problem}[theorem]{Problem}
\newcommand{\bpr}[2][nothing]{\ifthenelse{\equal{#1}{nothing}}
 {\begin{problem}} {\begin{problem}[#1]}\llabel{pr:#2}}

\newtheorem{corollary}[theorem]{Corollary}

\newcommand{\bpf}[1][Proof.]{\smallskip\noindent{\it #1} }
\newcommand{\qed}{\nolinebreak\mbox{\hspace{5 true pt}%
  \rule[-0.85 true pt]{3.9 true pt}{8.1 true pt}}}

\newcommand{\epf}{\qed \medskip}

\newcommand{\brm}{\smallskip\noindent{\bf Remark.} }

\newcommand{\OPdata}{Oleg Pikhurko\footnote{Partially supported by  the
National Science Foundation,  Grant 
DMS-0758057, and the Berkman Faculty Development Fund, CMU.}\\
Department of Mathematical Sciences\\
Carnegie Mellon University\\
Pittsburgh, PA 15213\\
Web: {\tt http://www.math.cmu.edu/\symbol{126}pikhurko}}

\begin{document}


\newcommand{\forb}{\mathrm{Forb}}
\newcommand{\Fi}{F}
\newcommand{\bydef}{: }
\newcommand{\weight}{\omega}
\newcommand{\btu}{{\bigtriangleup}}
\newcommand{\Id}{\mathrm{Id}}
\newcommand{\lmu}{\mu}
\newcommand{\ex}{\mathrm{ex}}
\newcommand{\Ex}{\mathrm{\C E\!\C X}}
\newcommand{\dd}{\hspace{1.0pt}\mathrm{d}\hspace{0.4pt}}
\newcommand{\esssup}{\mathrm{ess\hspace{1.1pt}sup}}
\newcommand{\Int}{_*}
\newcommand{\UInt}{^*}
\newcommand{\Lim}{\mathrm{LIM}}

\newcommand{\W}{\C W}
\newcommand{\WO}{\C W_{I}}
\newcommand{\A}[1]{A(#1)}
\newcommand{\CX}{{\C X}}
\newcommand{\const}[1]{\mathrm{Const}(#1)}
\newcommand{\continuum}{\mathfrak{c}}
\newcommand{\IR}{\I R}

\title{An Analytic Approach to Stability}
\author{\OPdata}
\maketitle

\begin{abstract}
 The stability method is very useful for obtaining exact
 solutions of 
many extremal graph problems. Its key step is to
establish the \emph{stability property} which, roughly speaking, states that
any two almost optimal graphs of the same order $n$ can be made isomorphic  by
changing $o(n^2)$ edges. 

Here we show how the recently developed
theory of graph limits can be used to give an analytic approach to
stability.  As an application, we present a new proof of the
Erd\H os--Simonovits Stability Theorem.

Also, we investigate various properties of the edit distance. In
particular, we show that the combinatorial and fractional versions
are within a constant factor from each other, thus answering a question of
Goldreich, Krivelevich, Newman, and Rozenberg.\end{abstract}

\section{Introduction}\llabel{intro}

The notion of the left convergence of graph sequences was
introduced by Borgs, Chayes, Lov\'asz, S\'os, and Vesztergombi (2003,
unpublished) and was developed
in~\cite{borgs+chayes+lovasz:10,BCLSV:06,BCLSV,BCLSV:08,diaconis+janson:08,elek+szegedy:08:arxiv,freedman+lovasz+schrijver:06,janson:08:arxiv,lovasz+szegedy:06,lovasz+szegedy:07:gafa,lovasz+szegedy:09:arxiv}
and other papers.  Benjamini and Schramm~\cite{benjamini+schramm:01}
introduced convergence for graphs of bounded maximum
degree. Tardos~\cite{tardos:trees} defined limits of
trees. Lov\'asz~\cite{lovasz:09} presents a nice survey of this area.

It is possible that graph limits will become a very powerful tool,
especially in extremal graph theory. The left limits are closely related to
the (Weak) Regularity Lemma, see Lov\'asz and
Szegedy~\cite{lovasz+szegedy:07:gafa}, which is a very important and
useful result. The algebraic characterization of Lov\'asz and
Szegedy~\cite[Theorem~2.2]{lovasz+szegedy:06} of possible limiting
subgraph densities seems to have a great potential.  Although these
developments are very recent, Razborov~\cite{razborov:07,razborov:08}
has already used graph limits to obtain a spectacular progress on
the long-standing Rademacher-Tur\'an problem.  Also, graph limits have
proved helpful for property and parameter testing, see Benjamini,
Schramm, and Shapira~\cite{benjamini+schramm+shapira:10}, Borgs et
al~\cite{BCLSSV:06}, Elek~\cite{elek:07:arxiv}, Lov\'asz and
Szegedy~\cite{lovasz+szegedy:08:arxiv}, and other.

Here is an example of how graph limits may be applied to extremal
graph problems.

Suppose that the convergence on graphs is encoded by a compact metric
space $(\CX,\delta)$ and a map that corresponds to each graph $G$ a point
$\A G$ of $\CX$ and \emph{respects graph isomorphism} (that is, $\A G=\A
H$
whenever $G\cong H$). Then we say that a sequence of graphs
$(G_n)_{n\in\I N}$ \emph{converges} if the sequence $(\A {G_n})_{n\in\I
  N}$ is Cauchy in the metric $\delta$. In this case, the \emph{limit} of
$(G_n)_{n\in\I N}$ is the (unique) limiting point of the sequence
$(\A {G_n})_{n\in\I N}$ in $(\CX,\delta)$, which exists since $(\CX,\delta)$ is
compact.

Suppose that we are given a \emph{graph parameter} $f$, that is, a
function on graphs that respects graph isomorphism, and
a \emph{graph property} $\C P$, that is,
a family of graphs closed under isomorphism. Let 
$\C P_n=\{G\in\C P\mid v(G)=n\}$
consist of all graphs in $\C P$ with $n$ vertices.
The corresponding \emph{extremal $(f,\C
P)$-problem} is to determine for each $n$
 \begin{eqnarray*}
 \ex_f(n,\C P)&=&\max\{f(G)\mid G\in \C P_n\},\\
 \Ex_f(n,\C P) &=& \{ G\in \C P_n\mid  f(G)=\ex_f(n,\C P)\},
 \end{eqnarray*}
 the maximum of $f(G)$ over all graphs from $\C P_n$ 
as well as the set of \emph{extremal graphs}, i.e.\ graphs that
achieve this maximum.
For example, if we let $h(G)$ be the maximum size of a \emph{homogeneous set} (a
clique or an independent set) in a graph
$G$, $f(G)=-h(G)/\log_2 v(G)$ be its scaled version, and 
$\C P$ be the family of all graphs, 
then we obtain the inverse problem for the diagonal Ramsey numbers.
Many
extremal graph problems can be represented this way.

Let us try to formulate some approximation (the ``limiting'' case) 
of the problem as
$n\to\infty$. We suggest the following definition. 
Let the \emph{limit set} $\Lim(f,\C P)$
consist of those $x\in \CX$ for which there is an infinite increasing  
sequence of indices $n_1<n_2<n_3<\dots$ and graphs 
$G_{n_i}\in\C P_{n_i}$ such that 
 \beq{AlmostExtr}
 \lim_{i\to\infty} (f(G_{n_i})-\ex_f(n_i,\C P))=0
 \eeq
 and the sequence $(G_{n_i})_{i\in \I N}$ converges to $x$, that is,
 $$
 \lim_{i\to\infty} \delta(\A {G_{n_i}},x)=0.
 $$
 Although we are ultimately interested in $\Ex_f(n,\C P)$, we do not
require that $G_{n_i}\in\Ex_f(n_i,\C P)$ here. One of
the reasons is that we often know $\ex_f(n,\C P)$ asymptotically but
not exactly, in
which case one can test if~\req{AlmostExtr} holds but not
the membership in $\Ex_f(n,\C P)$.

Now, we can try to study the set  $\Lim(f,\C P)$, 
which is independent of
$n$. If we succeed in completely describing it, 
then we might be able to discover some information about extremal
graphs. Indeed, if we select arbitrary extremal graphs
$G_n\in \Ex_f(n,\C F)$ for infinitely many $n$, then, by the compactness of $(\CX,\delta)$,
there always is a convergent
subsequence, whose limit belongs to $\Lim(f,\C
P)$. Suppose that this convergence implies 
some structural statement (in
purely graph theoretical terms)
that
necessarily occurs for infinitely many of the selected extremal
graphs. Then one can conclude that the statement fails only for
finitely many extremal graphs overall. 

One can call this approach the \emph{limit method}. It applies in
principle to very general settings. For example, the families $\C
P_n$ need not be related to each other for different $n$ nor the graph
parameter $f$ has to behave well with respects to taking limits: the
above definitions make perfect sense for arbitrary $f$ and $\C P$ (and
$\Lim(f,\C P)\not=\emptyset$ provided 
infinitely many of $\C P_n$'s are non-empty). Also, the
definition of the limit set may be modified to work with other extremal
problems, those which are indexed by a different parameter than the order of
a graph.

Since the limit method deals only with some approximation of the extremal
problem, one would hope to obtain only the
asymptotic of $\ex_f(n,\C P)$ at best. However, this approach might work
well together with the so-called \emph{stability method} that has proved very
useful in solving many extremal problems exactly (including the description
of $\Ex_f(n,\C P)$) for all large $n$.

The \emph{stability method} proceeds as follows.  Suppose that we know
the value of $\ex_f(n,\C P)$ asymptotically and that we have some set
$\C C_n$ believed to be exactly the set $\Ex_f(n,\C P)$ for all large
$n$. Assume that $\C C_n\subseteq \C P_n$ and $f$ is constant on $\C
C_n$. (Of course, these assumptions are necessary for $\C
C_n=\Ex_f(n,\C P)$ and, usually, they are easy to check.)  Given $\C
C_n$, we have to prove first that for any almost
extremal graph $G\in \C P_n$ (i.e.\ $G\in\C P_n$ satisfying
$f(G)=\ex_f(n,\C P)-o(1)$) there is $H\in\C C_n$ such that
$\hat\delta_1(G,H)=o(1)$, where
  \beq{hatdelta1}
 \hat\delta_1(G,H)=\frac{2}{n^2}\, \min\{|E(G)\bigtriangleup
\sigma(E(H))|\mid \mbox{bijective }\sigma:V(H)\to V(G)\}
 \eeq 
 is  the \emph{edit distance} between two
graphs of the same order $n$: it is $2/n^{2}$ times
the minimum number of adjacencies that one has to change in $G$ to
make it isomorphic to $H$. 
Next, pick an arbitrary $G\in\Ex_f(n,\C
P)$ for a sufficiently large $n$. 
By the above, we know that $G$ is close in the distance $\hat\delta_1$ 
to the graph property $\C
C_n$. In order to complete the proof, it is enough to argue that $G$
is necessarily in $\C C_n$. Here we can use various arguments, such as
applying ``local improvements'' to $G$ or arguing that every ``wrong'' adjacency
in $G$ bears too much penalty. Knowing all but $o(n^2)$ edges of $G$
greatly helps in this task; this is what makes this method so
successful. This approach was pioneered by
Simonovits~\cite{simonovits:68} in the late 1960s. It has been used to
obtain exact solutions for an impressive array of problems since then.

The term ``stability'' refers to the property that every almost extremal
graph has structure almost the same as some extremal graph. A
class of extremal problems for which this method seems to be particularly
suited is when there is only one pattern independent of $n$ 
for all almost extremal graphs. In order to state this property
formally, we have to define a  
version of edit distance for arbitrary pairs of graphs. 
Namely, the \emph{$\delta_1$-distance}, denoted by $\delta_1(G,H)$,
between graphs
$G$ and $H$ on vertex sets $\{x_1,\dots,x_m\}$ and $\{y_1,\dots,y_n\}$
respectively 
is the minimum over all non-negative $m\times n$-matrices
$A=(\alpha_{i,j})$ with row sums $1/m$ and column sums $1/n$ of
 \beq{GH}
 \delta_1(G,H,A)=\sum_{(i,j,g,h)\in \bigtriangleup}
 \alpha_{i,g}\alpha_{j,h},
 \eeq
 where $\bigtriangleup$ consists of all quadruples
$(i,j,g,h)\in [m]^2\times [n]^2$ 
such that exactly one of the following two relations
holds: either $\{x_i,x_j\}\in E(G)$ or $\{y_g,y_h\}\in
E(H)$. Informally speaking, we view $G$ and $H$ as uniformly vertex-weighted
graphs of total weight $1$ while 
$\alpha_{i,j}$ tells what fraction of
vertex $x_i$ is mapped into vertex $y_j$. It is not hard to show (see Section~\ref{delta1}) that
this defines a \emph{pre-metric} on the set of graphs, that is,
$\delta_1$ is  symmetric, non-negative and satisfies the
Triangle Inequality (but  may assume value zero 
on distinct graphs: e.g.\ $\delta_1(K_{m,m},K_{n,n})=0$ for any $m,n>0$).

Note that, for graphs $G_1$ and $G_2$ of
the same order, we trivially have
$\hat\delta_1(G_1,G_2)\ge \delta_1(G_1,G_2)$. This inequality is in general 
strict (see Arie
Matsliah's example presented in the technical report~\cite[Appendix
B]{goldreich+krivelevich+newman+rozenberg:08}
or Example~\ref{ex:diff} here). However,
we prove in Lemma~\ref{lm:relate} that
 \beq{relate}
  \hat \delta_1(G_1,G_2)\le 3\,\delta_1(G_1,G_2),
 \eeq
 answering in the affirmative an open question posed by Goldreich,
Krivelevich, Newman, and Rozenberg~\cite[Section
6]{goldreich+krivelevich+newman+rozenberg:08}
(see~\cite{goldreich+krivelevich+newman+rozenberg:09} for the journal version).

Now, let us say that
the extremal $(f,\C P)$-problem is \emph{stable} if for every $\e>0$
there are $\e'>0$ and $n_0$ such that for every $n_1,n_2\ge n_0$ and every
two graphs $G_1,G_2$ with $G_i\in \C P_{n_i}$ and $f(G_i)\ge 
\ex_f(n_i,\C P)-\e'$, for $i=1,2$, we necessarily have
$\delta_1(G_1,G_2)<\e$.
Theorem~\ref{th:stable} here gives an alternative
characterization of stable extremal problems. 
However, we postpone the
exact statement as well as the proof until Section~\ref{stab} after
we define graph limits in Section~\ref{graphons} and extend the
distance $\delta_1$ to them in Section~\ref{delta1}.

For example, our approach applies to the \emph{Tur\'an problem} that
asks for the maximum size of an $\C F$-free graph of order $n$.  This
is a central question of extremal graph theory that was introduced by
Tur\'an~\cite{turan:41}. Its scaled version can be represented in our
notation as $\ex_\rho(n,\forb(\C F))$, where
$\rho(G)=2e(G)/(v(G))^2$ denotes the \emph{edge density} of $G$ and
$\forb(\C F)$ consists of all $\C F$-free graphs. By applying our
Theorem~\ref{th:stable}, we obtain a new proof of the following
celebrated result in Section~\ref{es}.

\bth[The Erd\H os--Simonovits Stability
Theorem~\cite{erdos:67a,simonovits:68}]{ES}
For every
(possibly infinite) family  $\C F$
of non-empty graphs the extremal $(\rho,\forb(\C F))$-problem is stable.\end{theorem}

It is well-known that $\ex_\rho(n,\forb(\C F))=\frac{r-1}r
+o(1)$, where 
 \beq{r}
 r=\min\{\chi(F)\mid F\in\C F\}-1\ge 1,
 \eeq
 and the lower bound is given by \emph{the
Tur\'an graph} $T_r(n)\in\forb(\C F)$, the complete $r$-partite graph on $[n]$
with parts of size $\floor{n/r}$ or $\ceil{n/r}$. Thus, by~\req{relate},
Theorem~\ref{th:ES} can be reformulated in the more familiar form that 
for any
$\e>0$ there are $\e'>0$ and $n_0$ such that
every $\C F$-free graph with $n\ge n_0$ vertices and at least
$(\frac{r-1}r- \e') {n\choose 2}$ edges can be made isomorphic to $T_r(n)$ by
changing at most $\e {n\choose 2}$ edges.

Theorem~\ref{th:ES} was first applied by
Simonovits~\cite{simonovits:68} to determine the exact value of the
Tur\'an function $\ex(n,F)$ for various forbidden graphs $F$. This
theorem has a huge number of  applications. For example, Theorem~\ref{th:ES} turned up quite a few times in the author's research
alone: see the papers with
Jiang~\cite{jiang+pikhurko:09}, 
Lazebnik and
Woldar~\cite{lazebnik+pikhurko+woldar:07},
Loh and Sudakov~\cite{loh+sudakov+pikhurko}, 
Mubayi~\cite{mubayi+pikhurko:07}, 
Yilma~\cite{pikhurko+yilma}. Another proof of
Theorem~\ref{th:ES} was recently discovered by F\"
uredi~\cite{furedi:07:misc}.

\section{Graph Limits}\llabel{graphons}

\newcommand{\IB}{I_{\C B}}

Here we present the main definitions of ``dense'' graph limits. This notion of convergence (also called the
\emph{left convergence} in~\cite[Section~2.2]{BCLSV:08}) will be of main
interest for this paper. We refer the reader to
e.g.~\cite{BCLSV:08} for further details.

Until recently, the measure-theoretic methods were rare in discrete
mathematics (if compared with, for example, linear algebra or
topological tools). Bearing in mind a combinatorialist reader who does
not use real analysis in research, we decided to take an extra
care with measure theoretic concepts and to give references or
detailed explanations whenever feasible (even of some fairly standard
results). For example, the result of Lemma~\ref{lm:0} is stated
in~\cite[Page~5]{lovasz+szegedy:08:arxiv} without proof; here we
carefully fill in all missing details.
All analytical terms that
we do not define can be found in the book by
Folland~\cite{folland:ra}.

Let $\I R$ denote the set of reals and $I\subseteq \I R$ denote the closed unit interval
$[0,1]$. 
 For $Y\subseteq \IR^n$, let $\C L_Y=\{A\cap Y\mid A\in\C L\}$ denote
the restriction of the $\sigma$-algebra $\C L$ of Lebesgue measurable
subsets of $\IR^n$ to $Y$. If $Y\subseteq \IR^n$ is Lebesgue
measurable, then $\lmu_Y$ denotes the restriction of the Lebesgue
measure $\lmu$ to $\C L_Y$. Let $\C B_Y=\{A\cap Y\mid A\in\C B\}$ be
the restriction of the $\sigma$-algebra $\C B$ of Borel subsets of $\IR^n$ to $Y$.  When the set $Y$ is clear from the context, we write $\C
L$, $\lmu$, and $\C B$ for $\C L_Y$, $\lmu_Y$, and $\C B_Y$
respectively.
We say that some property holds \emph{almost
everywhere} (abbreviated as \emph{a.e.}) if the set of $x$ for which it
fails has Lebesgue measure 0. A measurable function 
is called \emph{simple} if it assumes only finitely many
values. 

A function $W:I^2\to \I
R$ is called \emph{symmetric} if $W(x,y)=W(y,x)$ for every $x,y\in I$. Let
$\W$ consist of all symmetric bounded measurable functions $W:(I^2,\C
L)\to(\I R,\C B)$. Following~\cite{BCLSV:08}, we call the elements of $\W$
\emph{graphons}. Let $\WO$ consist of those graphons $W\in\W$ such that $0\le
W(x,y)\le 1$ for every $x,y\in I$.  

A function $\phi:(I,\C L,\lmu)\to (I,\C L,\lmu)$ is called \emph{measure
preserving} if it is measurable and $\lmu(\phi^{-1}(A))=\lmu(A)$ for any
$A\in\C L_I$. Let $\Phi$ consist of all such functions. 
Note that $\phi\in \Phi$ may
be very far from being invertible as e.g.\ $\phi(x)=2x-\floor{2x}$ shows.
Let $\Phi_0$ consist of bijections
$\phi:I\to I$ such that both $\phi$ and $\phi^{-1}$ belong to $\Phi$.
Clearly, each of $\Phi$ and $\Phi_0$ is closed under taking compositions of functions.
For $\phi\in\Phi$ and $W\in \W$, let $W^\phi$ be defined by
$W^\phi(x,y)=W(\phi(x),\phi(y))$. It is easy to see that $W^\phi\in
\W$ and for any $\psi\in \Phi$, we have 
 \beq{circ}
 (W^{\phi})^\psi=W^{(\phi\circ \psi)}.
 \eeq

A few remarks are in order. It is standard (see e.g~\cite[Page
44]{folland:ra}) to consider the $\sigma$-algebra $\C B$ of
\emph{Borel} sets whenever (a subset of) $\IR^n$ is the range of a
function from some measure space. This has many
advantages: we can add or multiply such
functions~\cite[Proposition~2.6]{folland:ra}, take pointwise
limits~\cite[Proposition~2.7]{folland:ra}, etc, with the resulting
function being measurable. 
In particular, by \cite[Theorem~6.6]{folland:ra}, the vector space 
 \beq{L1}
 L^1:=L^1(I^2,\C L,\lmu)= \left\{\,\mbox{integrable $W:(I^2,\C
     L,\lmu)\to(\IR,\C B)$}\,\right\}/\sim,
 \eeq
 where we write $U\sim W$ iff $U=W$ a.e.,
is a Banach space with respect to the \emph{$\ell_1$-norm} 
 \beq{ell1}
 \|W\|_1=\int_{I\times I} |W(x,y)|\, \dd \lmu(x,y).
 \eeq
 On the other hand,
DiBenedetto~\cite[Section 14.1]{dibenedetto:ra}
demonstrates that the set of measurable functions from $(I,\C L)$ to $(I,\C L)$ is
not closed under taking pointwise 
limits (nor under multiplication, nor under addition, even if we take the
interval $[0,2]$ as the new range, as some easy modifications of his example can
show).  Note that, by definition, the set $\Phi$ consists of Lebesgue-to-Lebesgue measurable
functions (so that, e.g., for every $W\in\WO$ and $\phi\in \Phi$, we have
$W^\phi\in \WO$).

One can show that for any $W\in\W$ there is $U\in \W$ such that $W=U$ a.e.\ and
$U$ is measurable as a function  $(I^2,\C B)\to (\I R,\C B)$. (Indeed, by
writing the values of $W$ in base $2$, 
represent $W=\sum_{i\in\I Z} 2^i\, \I I_{X_i}$ as a linear
combination of the indicator functions of Lebesgue 
sets $X_i\in\C L_{I^2}$ and then
replace each $X_i$ by some Borel set
$Y_i\in \C B_{I^2}$ with $\lmu(X_i\bigtriangleup
Y_i)=0$.) This allows some flexibility in the definitions above. Still, in
order to eliminate any ambiguity, we decided to specify the corresponding
$\sigma$-algebras whenever the measurability of functions may matter.

Also, note that every graphon $W\in\W$, as a bounded
measurable function on the finite measure space $(I^2,\C L,\lmu)$, 
is integrable (see~\cite[Section~2.2]{folland:ra}), that is, $W\in
L^1$. 

Finally, let us remark that the standard definition of $L^1$ allows
functions to assume values $\pm\infty$. (This is convenient in the
statements of many theorems
of real analysis.) Since any integrable function assumes value
$\pm\infty$ on a set of measure $0$ and we identify a.e.\ equal
functions, we can restrict ourselves in~\req{L1} 
to functions with values in $\I
R$ only. 

For any integrable function $W:(I^2,\C L,\lmu)\to (\IR,\C B)$ (in particular,
for any graphon $W$), define its 
\emph{cut-norm} (also called the \emph{box-norm},
\emph{rectangle-norm}, etc) by
 \beq{cut}
 \|W\|_\Box=\sup_{S,T\in \C L_I} \left| \int_{S\times T} W(x,y)\, \dd
   \lmu(x,y)\right|.
 \eeq
 The \emph{cut-distance} $\delta_\Box(U,W)$ between $U,W\in \W$ is the infimum of
$\|U-W^\phi \|_\Box$ over all $\phi\in\Phi_0$. 
See~\cite[Lemma~3.5]{BCLSV:08} for other equivalent ways to
define this distance. For any $S,T\in\C L_I$, $\phi\in\Phi_0$, and an
integrable function $W:(I^2,\C L,\lmu)\to(\I
R,\C B)$,  
 we have 
 \beq{phi}
 \int_{S\times T} W(x,y)\,
\dd\lmu(x,y) = \int_{\phi^{-1}(S)\times \phi^{-1}(T)} W^\phi(x,y)\,
\dd\lmu(x,y),
 \eeq which is easiest to see from the definition of the
Lebesgue integral by approximating $W$ by simple
functions~\cite[Section~2.2]{folland:ra}. It follows that
$\|U-W^\phi\|_\Box=\|U^{\phi^{-1}}-W\|_\Box$ and that
$\delta_\Box$ is a pre-metric on $\WO$ (see the
argument leading to~\req{TI}). 

For a graphon $W\in\W$ we consider its equivalence class
 $$
 [W]=\{\,U\in \W\mid \delta_\Box(U,W)=0\,\}.
 $$
  Let
 \beq{X}
 \CX=\{\,[W]\mid W\in \WO\,\}
 \eeq
 consist of those equivalence classes that have a representative 
in $\WO$. We call elements of $\CX$  \emph{graph
limits}.  The pre-metric
$\delta_\Box$ induces a metric on $\CX$, which we still denote 
by the same symbol $\delta_\Box$.

Usually, it is more convenient to operate
with graphons, understanding equivalence classes implicitly.
But here we try to be as explicit as it is reasonably possible. 
Since  the words ``graph'' and ``limit'' are frequently used in this paper in various
contexts, we will use (in the absence of a better name) 
the term \emph{graphit} 
when referring to an equivalence
class $[W]$ with
$W\in\WO$. (One might view terms ``graphon'' and ``graphit'' as abbreviations of
``graph function'' and ``graph limit''.) 

For a graph $G$ on vertices $\{x_1,\dots,x_n\}$,
the corresponding element of
$\CX$ is $\A G=[W_G]$, where  $W_G\in \WO$  is defined by
 \beq{WG}
 W_G(x,y)=\left\{\begin{array}{ll} 1, & \mbox{if $(x,y)\in
       [\frac{k-1}n,\frac kn)\times [\frac{l-1}n,\frac ln)$ and
$\{x_k,x_l\}\in E(G)$,}\\
0,& \mbox{for all other $(x,y)\in I^2$,}
 \end{array}\right.
 \eeq
 that is, we encode the adjacency matrix of $G$ by a function $W_G\in\WO$.
  Clearly, the graphit $\A G$ does not depend on the labeling of $V(G)$ (while
the graphon $W_G$ does in general).

We have completely defined the metric space $(\CX,\delta_\Box)$ and the
special points $\A G$. This determines the promised convergence on
graphs. Let us give some brief pointers to the main properties of this
construction. 

Lov\'asz and Szegedy~\cite[Theorem~5.1]{lovasz+szegedy:07:gafa} proved
that the metric space $(\CX,\delta_\Box)$ is compact. Also, they
showed 
\cite[Theorem~2.2]{lovasz+szegedy:06} that the set $\{\,[W_G]\mid
\mbox{$G$ is a graph}\,\}$ is dense in $(\CX,\delta_\Box)$, that is,
every graphit $[W]$ with $W\in \WO$ is a limit of some sequence of graphs.

Any graph sequence $G_n$ with $e(G_n)=o(v(G_n)^2)$ as $n\to\infty$,
converges to the graphit $[\const0]$, where for $\alpha\in I$,
$\const{\alpha}\in\WO$ is the constant function that assumes the value
$\alpha$. This is why the phrase ``convergence of dense graphs'' is
often used.

The graphon $W_G$ can be viewed as a version of the adjacency matrix
of a graph $G$. However, a better informal interpretation of a general
graphon $W\in\WO$ is as a continuous version of the matrix that
encodes densities between parts of a (weak) regularity partition, see
\cite[Section~5]{lovasz+szegedy:07:gafa}. This also hints why,
although we start with $0/1$-valued functions $W_G$, we have to
allow general real-valued functions when we pass to limits. Having this data
for the graph, one can approximate, for example, the value of a
max-cut: for graphons the corresponding
computation is the supremum of the integral in~\req{cut} over
disjoint measurable $S,T\subseteq I$.

For graphs $F$ and $G$ the \emph{density} $t(F,G)$
of $F$ in $G$ is the probability that a random (not necessarily
injective) map $V(F)\to V(G)$
induces a homomorphism from $F$ into $G$.

As it turns out, the subgraph densities behave well with respect to
the $\delta_\Box$-distance. In combinatorial terms, this says, roughly speaking, that if
for two graphs $G$ and $H$ on $[n]$ we have
 \beq{SmallCuts}
 \big|\,e(G[A,B])-e(H[A,B])\,\big|=o(n^2),\qquad \mbox{for every $A,B\subseteq
   [n]$},
 \eeq
 then for every fixed
graph $F$ we have $|t(F,G)-t(F,H)|=o(1)$.  We refer the reader to
\cite[Lemma~4.1]{lovasz+szegedy:06} or
\cite[Theorems~2.3 and 3.7]{BCLSV:08} for the precise statements and
proofs. This may be viewed as a version of the Counting Lemma: if we
know the pairwise densities in a regularity partition
$V(G)=V_1\cup\dots\cup V_k$ of a graph $G$,
and generate the corresponding $k$-partite random graph $H$ on $V(G)$, 
then as $v(G)$ and $k$ tend to
infinity, with high probability~\req{SmallCuts} holds, and we can approximate 
subgraph densities in $G$ by those in $H$. 
 This greatly
motivates why the cut-norm is chosen to define the distance on
graphons. The role of $\phi$ in the definition of $\delta_\Box$ is, in
the discrete language, to overlay fractionally the vertex sets of two
graphs, cf~\req{GH} here and~\cite[Section~5.1]{BCLSV:08}.

It is natural to define the \emph{density} of a graph $F$ on $[k]$ in
a graphit $[W]$ by picking an arbitrary graph sequence $(G_n)$
convergent to $[W]$ and letting 
 \beq{tFWLimit}
 t(F,[W])=\lim_{n\to\infty} t(F,G_n).
 \eeq 
 This is well-defined and does not depend on the choice of
$(G_n)$. In fact, by writing $t(F,G_n)$ as a $k$-fold sum and
approximating it by a $k$-fold integral, one can show (see
\cite[Lemma~4.1]{lovasz+szegedy:06} or
\cite[Theorem~3.7.a]{BCLSV:08}) that
 \beq{tFW}
 t(F,[W])=\int_{I^k} \prod_{\{i,j\}\in E(F)} W(x_i,x_j)
 \,\dd\lmu(x_1,\dots,x_k).
 \eeq 
 Furthermore, neither of these definitions depends on the choice of $W\in [W]$, so we
can write $t(F,W)$ in place of $t(F,[W])$. Also, we have $t(F,G)=t(F,W_G)$.

More generally, in terms of graphons,
\cite[Lemma~4.1]{lovasz+szegedy:06} (see also \cite[Theorem~3.7.a]{BCLSV:08})
implies that the induced function
$t(F,{-}):(\CX,\delta_\Box)\to I$ is continuous for any $F$. Thus if
$(W_n)_{n\in\I N}$ 
is $\delta_\Box$-Cauchy, then the sequence 
$(t(F,W_n))_{n\in\I N}$ of reals is Cauchy for every fixed graph $F$.
The converse of
this also holds, by a result of Borgs et al
\cite[Theorem~3.7.b]{BCLSV:08}. Thus for $W,W_1,W_2,\dots \in \W_I$,
 \beq{subgraphs}
 \lim_{n\to\infty} \delta_\Box(W_n,W)=0\qquad\mbox{if and only if}\qquad 
 \forall \mbox{ graph}\ F\quad \lim_{n\to\infty} t(F,W_n)=t(F,W).
 \eeq
 It follows that
 each graphit $[W]$ is uniquely
determined by its ``moments function'' $t({-},W)$. An algebraic characterization of
all possible functions $t({-},W)$ realizable by some $W\in\WO$ is given
by Lov\'asz and Szegedy~\cite[Theorem~2.2]{lovasz+szegedy:06}.

Let us also say a few words about graph limits and property testing.
(See Goldreich, Goldwasser, and Ron~\cite{goldreich+goldwasser+ron:98} for a precise definition of
property testing and several fundamental results.) In the most 
restrictive
sense (the \emph{oblivious} or \emph{order independent testing}), 
we have a (very big)
unknown graph $G$ and are told the subgraph $G[X]$ induced by a
random $m$-set $X$ of vertices, where $m$ is a fixed number. 
It is known that with probability at least $1-\e$ we
have $\delta_\Box(W_{G[X]},W_G)\le \e$, provided $m\ge m_0(\e)$
(see \cite[Theorem~2.5]{lovasz+szegedy:06} or 
\cite[Theorem~3.7]{BCLSV:08}). This means that we can learn a good
$\delta_\Box$-approximation to the graph $G$. The objective of
\emph{property testing} is to approximate with high probability how
far 
$G$ is from a given property $\C P$, but the \emph{edit} distance
$\hat\delta_1$ is to be used here. Graphons seem to 
provide very convenient tools and language for dealing with
this problem
(which essentially amounts to relating the $\hat\delta_1$ and
$\delta_\Box$ distances from an arbitrary graph to the given
property), see~\cite{BCLSSV:06,lovasz+szegedy:08:arxiv}.

\section{Extending the $\delta_1$-Distance to Graph Limits}\llabel{delta1}

Here we show how to extend the distance $\delta_1$ from graphs to
graphits. This definition is standard but it seems that no formal proofs
of some of its properties have appeared in the literature. Therefore we give
careful proofs of all claims (or references to them).
The author thanks L\'aszl\'o Lov\'asz for pointing out that Lemma~\ref{lm:0} 
can be deduced from the results
in~\cite{borgs+chayes+lovasz:10,lovasz+szegedy:08:arxiv}, which is the proof presented
here.

Here is the definition of $\delta_1$ for graphits. First, we define $\delta_1$
on $\W$, the set of graphons. For $U,W\in\W$, let 
 \beq{delta1:0}
 \delta_1(U,W)=\inf\left\{\|U-W^\phi\|_1\mid \phi\in\Phi_0\right\},
 \eeq
 where $\|U-W^\phi \|_1$ is the standard $\ell_1$-norm of $U-W^\phi$ 
as defined
 by~\req{ell1}.

Clearly, $\delta_1$ is non-negative. It is symmetric by~\req{phi}. 
Also, $\delta_1$ satisfies the Triangle
Inequality. Indeed, for every $U,V,W\in \W$ and $\e>0$ we can choose
$\phi,\psi\in \Phi_0$
such that $\|U^\phi-V\|_1\le \delta_1(U,V)+\e$ and $\|V-W^\psi\|_1\le
\delta_1(V,W) +\e$. Now, by the Triangle Inequality for the $\ell_1$-norm,
 \begin{eqnarray}
 \delta_1(U,W) &\le& \| U-(W^\psi)^{\phi^{-1}}\|_1\ =\
 \|U^\phi-W^\psi\|_1 \nonumber\\
 &\le& \|U^\phi-V\|_1 +
\|V-W^\psi\|_1\ \le\  \delta_1(U,V) + \delta_1(V,W) +2\e.\llabel{eq:TI}
 \end{eqnarray}
 Since $\e>0$ was arbitrary, the claim follows. Hence, $\delta_1$ is a
pre-metric on $\WO$.

We will present an equivalent definition of
$\delta_1$ in Lemma~\ref{lm:delta1} and will conclude in Corollary~\ref{cr:delta1} that
$\delta_1$ gives a metric on $\CX$. Let us state
a few auxiliary or related results first.

\blm{0ae} Let an integrable
$W:(I^2,\C L,\lmu)\to(\IR,\C B)$ satisfy $\|W\|_\Box=0$. Then
$W=0$ a.e. In particular, for any $U,W\in\W$, $\|U-W\|_\Box=0$ implies that
$\|U-W\|_1=0$.  
\end{lemma}
 \bpf Let $Z$ be the \emph{Lebesgue set} of the function $W$, which
can be defined
as the set of those $(x,y)$ in the interior of $I^2$ such that
 \beq{LDT}
 \lim_{c\to 0\atop c>0}\ \frac{1}{\lmu(R_{x,y,c})}\, \int_{(x',y')\in R_{x,y,c}}
\Big|\,W(x',y')-W(x,y)\,\Big|\,\dd\lmu(x',y') = 0,
 \eeq
 where $R_{x,y,c}$ is the open rectangle $(x-c,x+c)\times (y-c,y+c)$.

The Lebesgue Differentiation Theorem (\cite[Theorem~3.21]{folland:ra})
implies that $\lmu(Z)=1$. If $W(x,y)\not=0$ for some $(x,y)\in
Z$, then by~\req{LDT} there is $c>0$ such that 
 $$
 \left|W(x,y)-\frac{1}{4c^2} \int_{R_{x,y,c}}
W\,\right|<\frac{|W(x,y)|}2.
 $$ 
 Thus $\|W\|_\Box\ge |\int_{R_{x,y,c}}
W\, |\ge 2c^2|W(x,y)|>0$, a contradiction. Thus $W=0$ a.e.\epf 

A function $U:I^2\to \I R$ is
called an \emph{interval step function} if there is a partition 
$I=I_1\cup\dots\cup I_k$ into finitely many intervals such that $U$ is
constant on each \emph{rectangle} $I_i\times I_j$. Any interval step function is a simple function.
Of course, such $U$ is necessarily measurable, even in
the strongest sense as a function from $(I^2,\C B)$ to $(\I R,2^{\I R})$.

\blm{approx} For any $\e>0$ and any integrable
function $W:(I^2,\C L,\lmu)\to (\IR,\C B)$ there is an interval
step function $U$ such that $\|W-U\|_1<\e$. Moreover, if $W\in\WO$, 
then we can also require that $U\in \WO$.
\end{lemma}
 \bpf
 The first part of the lemma follows from
\cite[Theorem~2.41]{folland:ra} (see
also~\cite[Lemma~3.2]{BCLSV:08}). Let us establish the second
part. Let $W\in\WO$ and $U_0$ be the interval
step function with $\|W-U_0\|_1<\e$, given by the first part. Let 
$U_1(x,y)=g(U_0(x,y))$, where $g(z)=\max(0,\min(1,z))$ maps $z\in
\I 
R$ to the nearest point from $I$. Since for every $z'\in I$ and $z\in
\I R$ we have $|g(z)-z'|\le |z-z'|$, we conclude that $\|
U_1-W\|_1\le\|U_0-W\|_1\le\e$. Finally, we take
$U(x,y)=(U_1(x,y)+U_1(y,x))/2$. Then the new interval step function
$U$ belongs to $\WO$. Also, in view of
inequality $|a-c|+|b-c|\ge 2\,|\frac{a+b}2-c|$ valid for any $a,b,c\in\I R$, 
we have
$\|W-U\|_1\le\|W-U_1\|_1\le \e$, as desired.\epf

\brm This approximation reminds the one given by the Weak Regularity
Lemma of Frieze and Kannan~\cite{frieze+kannan:99} (see
also~\cite[Section~2]{lovasz+szegedy:07:gafa}) with respect to the
cut-norm, except we cannot bound the number of parts in
Lemma~\ref{lm:approx} in terms of $\e$ only. This is an important
distinction between the cut-norm and the $\ell_1$-norm, giving another
motivation for taking $\delta_\Box$ as the distance between
graphons. This allows one to construct a finite $\e$-net for the metric
space $(\C
X,\delta_\Box)$. Namely, let $n=n(\e)$ be large and take all interval steps
functions with steps $[\frac in,\frac{i+1}n)$ that assume
values in
$\{\frac1n,\dots,\frac nn\}$; there are at most
$n^{n^2}<\infty$ such functions. Thus $(\C X,\delta_\Box)$ is
totally bounded,
which is one of the ingredients needed for compactness.
See~\cite[Theorem~5.1]{lovasz+szegedy:07:gafa} for more details.\medskip

\blm{Phi0} Let $X,Y\in \C L_I$ have measure $1$ and let $\psi$ be a
bijection from $X$ onto $Y$ such that for any interval $J\subseteq I$ the sets
$\psi(J\cap X)$ and
$\psi^{-1}(J\cap Y)$ are Lebesgue measurable with $\lmu(\psi(J\cap X))=\lmu(\psi^{-1}(J\cap Y))=\lmu(J)$. Then 
there is $\phi\in\Phi_0$ such that $\phi=\psi$ a.e.\end{lemma}
 \bpf Suppose first that $|I\setminus X|=|I\setminus Y|=\continuum$,
that is, the cardinality of both $I\setminus X$ and
$I\setminus Y$ is continuum. 
Let $\phi$ be an arbitrary bijection between $I\setminus X$ and
$I\setminus Y$ while $\phi(x)=\psi(x)$ if $x\in X$. Then $\phi=\psi$
a.e. Also, for any interval $J\subseteq I$, the pre-image $\phi^{-1}(J)$
differs from $\psi^{-1}(J\cap Y)\in \C L$ on a set of measure $0$, so
it is Lebesgue measurable of measure $\lmu(J)$. Since $\C B$ is
generated by intervals as a $\sigma$-algebra
(\cite[Theorem~1.6]{folland:ra}), it follows (e.g.\ by application
of the uniqueness claim of~\cite[Theorem 1.14]{folland:ra}) that $\phi$ is a measure
preserving function from $(I,\C L)$ to $(I,\C B)$. But a subset of $I$
is Lebesgue measurable set if and only if it can be sandwiched between
two Borel sets of the same measure
(\cite[Theorem~1.19]{folland:ra}). This easily implies that $\phi$ is
a measure preserving map from $(I,\C L)$ to $(I,\C L)$, that is, $\phi\in
\Phi$. Likewise, $\phi^{-1}\in \Phi$, giving $\phi\in \Phi_0$ as required.

Finally, suppose that, for example, $|I\setminus X|<\continuum$. 
Let $C\subseteq I$ be \emph{the Cantor
set}, which has measure 0 and cardinality continuum
\cite[Proposition~1.22]{folland:ra}. Let $X'=X\setminus C$ and
$Y'=Y\setminus \psi(X\cap C)$. Then $\psi$ maps $X'$ bijectively onto
$Y'$. Also, $\lmu(\psi(X\cap C))=0$. Indeed, for every $\e>0$, we can
find a set $J\supseteq C$ which is the union of finitely many
intervals of total length at most $\e$ that covers $C$. By the
assumption of the lemma, $\psi(X\cap J)$ has measure at most
$\e$. Since $\e>0$ was arbitrary, $\lmu(\psi(X\cap C))=0$. Thus
$\lmu(X\setminus X')=\lmu(Y\setminus Y')=0$ and the restriction
$\psi|_{X'}$  satisfies
the assumptions of the lemma. Since $|I\setminus X'|=|I\setminus
Y'|=\continuum$, we already know how to find the required $\phi\in\Phi_0$
for $\psi|_{X'}$. The very same function $\phi$ works for $\psi$ as well.\epf

Let us call a point $x$ lying inside a Lebesgue set $A\subseteq \I R$ a \emph{density point of $A$} if 
 $$
 \lim_{c\to0\atop c>0}\
\frac{\lmu(A\cap (x-c,x+c))}{2c}=1,
 $$ 
 or equivalently, if 
$x$ belongs to the
Lebesgue set (as defined by the $1$-dimensional version of~\req{LDT})
of the characteristic function $\I I_A:\I R\to\{0,1\}$ of $A$.
Again, Theorem~3.21 in~\cite{folland:ra} implies that almost every
point of $A$ is a density point.

The arithmetic operations and the linear order on $I=[0,1]$ play no
role in the definition of graphons;
see~\cite[Section~2.1]{borgs+chayes+lovasz:10} for a more general
point of view. The following simple lemma suffices for our purposes.

\blm{new} For every partition of $I=A_1\cup\dots\cup A_k$ into Lebesgue
measurable sets $A_i$ there are a partition $I=I_1\cup \dots\cup I_k$
into intervals and $\psi\in \Phi_0$ such that
$\lmu(\psi(A_i)\bigtriangleup I_i)=0$ for each $i\in [k]$.\end{lemma}

\bpf It is enough to prove the case $k=2$ with the general claim
following by a simple induction on $k$. 

Let $a_1=0$, $a_2=\lmu(A_1)$, $I_1=[0,a_2]$, and $I_2=I\setminus
I_1$. Assume that $0<a_2<1$ (for otherwise any $\psi\in \Phi_0$
works). 

Let $i=1$ or $2$.  Let
$X_i\subseteq A_i$ be the set of density points of 
$A_{i}$. For $x\in X_{i}$ let 
 $$
 \psi(x)=a_{i}+\lmu(A_{i}\cap [0,x]).
 $$
 Then
$\psi(X_i)$ lies in the interior of $I_{i}$. Indeed, if,
for example, $\psi(x)=a_i$, then $\lmu(A_{i}\cap
(-\infty,x))=0$, so $x$ cannot be a density point for
$A_{i}$. Likewise, if $y\in X_{i}\setminus\{x\}$ is another density point of
$A_{i}$, then $\psi(y)\not=\psi(x)$. Let $Y_{i}=\psi(X_{i})$.
The pre-image under 
$\psi$ of any open interval $J=(a_{i},a_{i}+b)\subseteq I_{i}$ is the
intersection of the interval $(0,c)$ with $X_{i}$, where
 $$
 c=\sup\{x\in I\mid \lmu(A_{i}\cap [0,x])< b\}= \sup\{x\in I\mid \lmu(X_{i}\cap [0,x])< b\}.
 $$
 Since $b\le\lmu(X_{i})$ and the measure $\lmu$ is continuous from below
(\cite[Theorem~1.8.c]{folland:ra}), we conclude that
$\lmu(\psi^{-1}(J))=b=\lmu(J)$. Also, for any open interval
$J=(b,c)\subseteq I$, the image under $\psi$
of $X_i\cap J$ is  $Y_{i}\cap J_{i}$, where 
 $$
 J_{i}=(a_{i}+\lmu(A_{i}\cap [0,b]),a_{i}+\lmu(A_{i}\cap [0,c]))
 $$
 is a subinterval of
$I_{i}$ with $\lmu(J_{i})=\lmu(J\cap
X_{i})$. 

Let $X=X_1\cup X_2$ and $Y=Y_1\cup Y_2$. 
It routinely follows that all assumptions of  Lemma~\ref{lm:Phi0} with
respect to the bijection $\psi:X\to Y$ are satisfied. The element
$\phi\in \Phi_0$ returned by Lemma~\ref{lm:Phi0} has the required properties.\epf

\blm{inverse} For every interval step function $U\in\W$ and $\phi\in\Phi$,
there is $\psi\in\Phi_0$ such that $(U^\phi)^\psi=U$ a.e.\end{lemma}
 \bpf Let $I=I_1\cup\dots\cup I_k$ be a partition into intervals such
that $U$ is constant on each rectangle $I_i\times I_j$. For
$i,j\in[k]$, let $\alpha_{i,j}=\lmu(A_{i,j})$, where $A_{i,j}=I_j\cap \phi^{-1}(I_i)$. Since
$\phi$ is measure preserving, $\sum_{j=1}^k \alpha_{i,j}=\lmu(I_i)$
for every $i\in [k]$. Partition the interval 
$I_i=I_{i,1}\cup\dots\cup I_{i,k}$ into
intervals of lengths respectively
$\alpha_{i,1},\dots,\alpha_{i,k}$. By Lemma~\ref{lm:new} find $\eta\in
\Phi_0$ such that $\lmu(\eta(A_{i,j})\bigtriangleup I_{i,j})=0$.  The element
$\psi=\eta^{-1}\in \Phi_0$ has the required properties
by~\req{circ} because for a.e.\ $x\in I_{i,j}$ we have $\psi(x)\in A_{i,j}$
and $\phi(\psi(x))\in I_i$.\epf

Lemmas~\ref{lm:approx} and \ref{lm:inverse} easily imply the following result.

\begin{corollary}\llabel{cr:inverse} For any $U,W\in\W$ and $\phi\in
\Phi$, we have $\delta_1(U,W)=\delta_1(U^\phi,W)$.\qed\end{corollary}

\bth{tfae} For $U,W\in \W$, the following are equivalent.

(a) For every graph $F$, we have $t(F,U)=t(F,W)$.

(b) $\delta_\Box(U,W)=0$.

(c) There are $\phi,\psi\in \Phi$ such that $U^\phi=W^\psi$
a.e.\end{theorem}
 \bpf The equivalence of (a) and (b) follows from~\req{subgraphs} (i.e.\ from
 \cite[Lemma~4.1]{lovasz+szegedy:06} and \cite[Theorem~3.7]{BCLSV:08}).
 The equivalence of (a) and (c) is
proved by Borgs, Chayes, and Lov\'asz~\cite[Corollary~2.2]{borgs+chayes+lovasz:10}.\epf

\blm{delta1} For any $U,W\in\W_I$, we have
 \beq{delta1}
 \delta_1(U,W)=\inf_{\phi,\psi\in \Phi} \|U^\phi-W^\psi\|_1.
 \eeq
 \end{lemma}
 \bpf Since $\Phi_0$ is a subset of $\Phi$ and $\Phi_0$ 
contains the identity
function $\Id:I\to I$, the ``$\ge$''-inequality in~\req{delta1} easily
follows. Let us show the converse.

Let $U,W\in\W$ and $\e>0$. By Lemma~\ref{lm:approx}
we can
find interval step functions $U_0$ and $W_0$ lying within $\e$ from respectively
$U$ and $W$ in the $\ell_1$-norm. For any $\phi,\psi\in \Phi$, we have by~\req{phi}
 $$
 \|U^\phi-W^\psi\|_1 \ge \|U_0^\phi - W_0^\psi\|_1-
\|U^\phi - U_0^\phi\|_1 - \|W^\psi-W_0^\psi\|_1\ge \|U_0^\phi - W_0^\psi\|_1
- 2\e.
 $$
 Likewise, $\| U - W^\phi\|_1\le \|U_0-W_0^\phi\|_1+2\e$. Since $\e>0$
was arbitrary, it is enough to prove~\req{delta1} on the additional
assumption that $U$ and $W$ are interval step functions.

Again, let $\e>0$. Let $\phi,\psi\in\Phi$ be such that
$\|U^\phi-W^\psi\|_1-\e$ is at most the right-hand side of~\req{delta1}.
By Lemma~\ref{lm:inverse}
choose $\eta\in\Phi_0$ such that $(W^\psi)^\eta=W$ a.e. Then, by~\req{circ},
 \beq{temp1}
 \|U^\phi-W^\psi\|_1 = \|(U^\phi)^\eta-(W^\psi)^\eta\|_1 = \| U^{(\phi\circ \eta)} -
W\|_1.
 \eeq
 Again, by Lemma~\ref{lm:inverse} applied to $U$ and $\phi\circ
\eta\in\Phi$, find $\nu\in\Phi_0$ such that $(U^{(\phi\circ
\eta)})^\nu=U$ a.e. From~\req{temp1} we conclude that
$\|U^\phi-W^\psi\|_1= \|U -W ^\nu\|_1$, which is at least the right-hand side of~\req{delta1:0}. Since $\e$ was arbitrary, the lemma follows.\epf

\blm{ext} For any two graphs $G$ and $H$, the  $\delta_1$-distance
$\delta_1(G,H)$ defined by~\req{GH} 
is equal to $\delta_1(W_G,W_H)$, where $W_G$ and $W_H$
are defined
by~\req{WG}.\end{lemma}
 \bpf Let $V(G)=\{x_1,\dots,x_m\}$ and $V(H)=\{y_1,\dots,y_n\}$. 
For $\phi\in \Phi_0$, $\|W_G^\phi-W_H\|_1$ equals to the expression
in~\req{GH} with $\alpha_{i,j}=\lmu(I_i\cap \phi^{-1}(J_j))$,
$I_i=(\frac{i-1}m,\frac im)$ and $J_j = (\frac{j-1}n,\frac jn)$. 
Conversely, given
numbers $\alpha_{i,j}$ such the matrix $(\alpha_{i,j})_{i,j\in
[n]}$ has row sums $1/m$ and column sums $1/n$, 
one can easily construct $\phi\in\Phi_0$ giving these $\alpha_{i,j}$ as
above.\epf

\blm{0} Let $U,W\in\W$ satisfy $\delta_\Box(U,W)=0$. Then
$\delta_1(U,W)=0$.\end{lemma} 

\bpf By Theorem~\ref{th:tfae}, there are $\phi,\psi\in\Phi$ such that
$U^\phi=W^\psi$ a.e. The claim follows by using the equivalent
definition of $\delta_1$ from Lemma~\ref{lm:delta1}.\epf


\begin{corollary}\llabel{cr:delta1} The function $\delta_1$ 
induces a metric on the set $\CX$ of graphits, extending the
$\delta_1$-distance from graphs.\qed\end{corollary}

\brm Let us point out that the convergence with respect to the cut-distance does not generally imply the convergence with
respect to $\delta_1$. For example, the infinite sequence of random graphs
$G_n\in G_{n,1/2}$ converges in the $\delta_\Box$-distance 
with probability 1 to
the graphit $[\const{1/2}]$ by~\cite[Corollary~2.6]{lovasz+szegedy:06}
while no graph sequence whatsoever can converge in the $\delta_1$-distance to
$[\const{1/2}]$ by Theorem~\ref{th:dichotomy} here.\medskip

\section{Comparing the Discrete and Fractional $\delta_1$-Distances}

Clearly, for graphs $G$ and $H$ of the same order we have
$\hat\delta_1(G,H)\ge \delta_1(G,H)$, where $\hat\delta_1$ is defined
by~\req{hatdelta1}.  The distances $\hat\delta_1$ and $\delta_1$ do not
coincide in general as Example~\ref{ex:diff} demonstrates. Independently, 
Arie Matsliah (see~\cite[Appendix
B]{goldreich+krivelevich+newman+rozenberg:08}) presented another construction that
achieves ratio $6/5$. Although our ratio is smaller (only $11/10$), the
ideas behind our construction are different from those 
of Matsliah and might be useful in the quest for better ratios.  
Hence, we decided to keep this example in the paper.

\begin{example}\llabel{ex:diff} There are graphs  $G$ and $H$ such that $v(G)=v(H)$ but
 $$
 \hat\delta_1(G,H)\ge \frac{11}{10}\, \delta_1(G,H)>0.
 $$
\end{example}
 \bpf  Fix an integer $n\ge 24$. Pick disjoint sets
$X=\{x_1,\dots,x_4\}$,
$M=M_1\cup\dots\cup M_4$, and $N=N_1\cup\dots\cup N_5$ 
with each $M_i$ having $4$ elements and each
$N_i$ having $n$ elements.

Let $V(G)=V(H)=N\cup M\cup X$. It will be the case that $N\cup M$
spans the same subgraph in both $G$ and $H$. Namely, $N$ spans the
complete graph while, for $i\in [4]$, we put the complete bipartite
graph between $M_i$ and $\cup_{j=1}^i N_j$. These are all edges
inside $M\cup N$. 

Fix another partition $M=L_1\cup\dots\cup L_4$ such that
each $L_i$ has $4$ elements and $|L_i\cap M_i|=|L_{i+1}\cap
M_{i}|=2$ for $i\in [4]$, where we agree that $L_5=L_1$.

In $G$, the edges incident to $X$ are as
follows: $\{x_i,x_j\}$ for $1\le i<j\le 4$ with $j-i$ even 
plus all pairs $\{x_i,y\}$
for $i\in[4]$ and $y\in M_i$. In $H$, the edges incident to $X$ are
as follows: $\{x_i,x_j\}$ for $1\le i<j\le 4$ with $j-i$ odd 
plus all pairs $\{x_i,y\}$ for $i\in [4]$ and $y\in L_i$. 

We have 
 \beq{optGG'}
 |E(G)\bigtriangleup E(H)|=  \sum_{i=1}^4
|M_i\bigtriangleup L_i|+ {|X|\choose 2} =22.
 \eeq
 
Let us show that this is smallest
possible. Pick an optimal bijection $\sigma:V(G)\to V(H)$. 
In each of $G$ and $H$, every vertex in $N$ has degree at
least $5n-1$ while any vertex in $M\cup X$ has degree at most
$4n+1$. Hence, if $\sigma$ does
not preserve $N$, then the number of discrepancies will be at least
$(5n-1)-(4n+1)\ge 22$. So, assume that $\sigma(N)=N$.
Likewise, we have $\sigma(M_i)= M_i$, for otherwise the number
of discrepancies (between $M$ and $N$) is at least $n>22$. Finally,
consider the action of $\sigma$ on $X$. For every $x,y\in X$, their
neighborhoods in $M$ with respect to $G$ and $H$ differ by at least
$4$. If $\sigma$ does not map some $x_i$ into $\{x_i,x_{i+1}\}$, where
$x_5=x_1$, then the neighborhoods $N_G(x_i)$ and $N_H(\sigma(x_i))$ in $M$
are disjoint and this vertex alone creates at least $8$
discrepancies. Moreover, since $X$ spans $2$ and $4$ edges in $G$ and $H$
respectively, the total number of discrepancies is at least $8+3\times
4+2=22$ and we cannot improve~\req{optGG'}. Thus let us assume that 
$\sigma(x_i)\in \{x_i,x_{i+1}\}$ for every $i\in
[4]$. This implies that either $\sigma$ is constant on $X$ or shifts
indices by $1$. In either case, this gives the same bound as
in~\req{optGG'}. 

Hence, $\hat \delta_1(G,H)\ge \frac{2\cdot 22}{(5n+20)^2}$. Let us establish
an upper bound on $\delta_1(G,H)$ now.

Let $G[2]$ be the \emph{2-fold blow-up} of $G$, where each vertex $x$
is replaced by two vertices $x',x''$ and each edge $\{x,y\}$ by the
complete bipartite graph with parts $\{x',x''\}$ and $\{y',y''\}$. For
$Y\subseteq V(G)$, let
$Y[2]=\{y',y''\mid y\in Y\}$.
Consider the following bijection $\sigma$ between the vertex sets of
$G[2]$ and $H[2]$. It is the identity bijection on $M[2]\cup
N[2]$.
For $i\in [4]$, let $\sigma(x_i')=x_i'$ and $\sigma(x_i'')=x_{i+1}''$. 
Easy checking shows that $\sigma$, when restricted
to $X[2]$, mismatches only $16$
adjacencies (versus $4\times {4\choose 2}=24$ if $\sigma$
were the identity). The number of discrepancies between $X[2]$ 
and $M[2]$ is
$4\times 16$. We have
 $$
 \delta_1(G,H)\le \hat\delta_1(G[2],H[2])\le \frac2{4(5n+20)^2}(4\times 16 
+16)\le \frac{10}{11}\, \hat\delta_1(G,H).\qed 
 $$

\blm{relate} For any two graphs $G$ and $H$ on the same vertex set 
$[n]$, we have
 $$
 \hat\delta_1(G,H)\le 3\delta_1(G,H).
 $$
 \end{lemma}
 \bpf If $G\cong H$, then
$\delta_1(G,H)=\hat\delta_1(G,H)=0$, 
so assume $G\not\cong H$. Let $\ell=n^2 \hat\delta_1(G,H)/2$ be the smallest
number of adjacencies we have to change in $G$ to make it isomorphic to~$H$.

Let $A=(\alpha_{i,j})_{i,j\in [n]}$ be an optimal overlay 
matrix as in~\req{GH}, where we assume $x_i=i$ and $y_j=j$. (Thus 
$nA$ is doubly-stochastic.) 

Although $nA$ can be represented as a convex combination of permutation
matrices by Birkhoff's theorem~\cite{birkhoff:46}, 
we find it more convenient to work with an
approximation where all coefficients are equal. (Thus some permutation
matrices may be repeated more than once.) Such an approximation is easy to
find as follows.

Pick a large $m>m_0(A)$. Inductively on $i$, we construct permutation matrices
$P_i$ as follows. Suppose that $i\ge 0$ and we have already found
$P_1,\dots,P_i$ such that $P'=P_1+\dots+P_i\le mnA$ (where matrix inequalities
are meant component-wise).  If there is a permutation matrix $P_{i+1}$ such
that $P'+P_{i+1}\le mnA$, take it and repeat the step.

Suppose that no such $P_{i+1}$ exists. Let $B=(\beta_{f,g})_{f,g\in
[n]}=mnA-P'$. This is a non-negative matrix with row/column sums $m-i$. By
Hall's Marriage theorem~\cite{hall:35}, 
there is a set $R\subseteq [n]$ of $r$ rows
and 
a set $S\subseteq [n]$ of $n-r+1$
columns such that each entry of the $R\times S$-submatrix of $B$ 
is less than 1. Hence,
 \begin{eqnarray*}
 (m-i)r&=&\sum_{f\in R}\sum_{g=1}^n \beta_{f,g}\ =\ \sum_{f\in
   R}\sum_{g\in S} \beta_{f,g} + \sum_{f\in R} \sum_{g\in
   [n]\setminus S}\beta_{f,g}\\
 &\le& r(n+1-r) + (m-i)(n-(n-r+1)),
 \end{eqnarray*}
 and therefore $m-i\le r(n+1-r)\le (n+1)^2/4$.  Let $P_{i+1},\dots,P_m$ be
arbitrary permutation matrices and  $P=\frac1{mn}(P_1+\dots+P_m)$.
It follows that
 $$
 \|A-P\|_\infty\le 2\times \frac{(n+1)^2}{4mn}=
\frac{(n+1)^2}{2mn}.
 $$
 Since $m$ is arbitrarily large, in order to prove the
lemma it is enough to show that 
 \beq{aim0}
 \hat\delta_1(G,H)\le 3\delta_1(G,H,P),
 \eeq
 where $\delta_1(G,H,P)$ is defined by~\req{GH}.

Let $\sigma_1,\dots,\sigma_m:[n]\to [n]$ be the permutations encoded by
$P_1,\dots,P_m$ respectively. 
As it was defined after~\req{GH}, $\btu$ is the set of all 
quadruples $(x,y,x',y')\in [n]^4$ such
that exactly one of the relations $\{x,y\}\in E(G)$ and $\{x',y'\}\in E(H)$
holds. Note that we allow $x=y$ or $x'=y'$ but both equalities cannot hold
simultaneously by the definition of $\btu$.

For $(i,j)\in [m]^2$, let $\btu(i,j)$ consist
of $(x,y)\in [n]^2$ 
such that $(x,y,\sigma_i(x),\sigma_j(y))\in \btu$. For
$(x,y,x',y')\in\btu$, let $I(x,y,x',y')$ consist of all pairs
$(i,j)\in [m]^2$ such that
$\sigma_i(x)=x'$ and $\sigma_j(y)=y'$.
 Also, for $X\subseteq [m]$, define
 $$
 S_X=\sum_{i,j\in X: i<j}  |\btu(i,j)|.
 $$
We have
 \begin{eqnarray}
 \delta_1(G,H,P)&=&\sum_{(x,y,x',y')\in \btu} P_{x,x'}P_{y,y'}\nonumber\\
  &=& \sum_{(x,y,x',y')\in \btu} \left(\sum_{i:\sigma_i(x)=x'}
   \frac1{mn}\right)\left(\sum_{j:\sigma_j(y)=y'}
   \frac1{mn}\right)\nonumber\\
 & =&  \frac1{m^2n^2}\sum_{(x,y,x',y')\in \btu} |I(x,y,x',y')|\nonumber\\
 &=&  \frac1{m^2n^2}\sum_{i,j\in [m]} |\btu(i,j)|
 \ =\ \frac{2S_{[m]} +\sum_{i=1}^m|\btu(i,i)|}{m^2n^2}.\llabel{eq:delta1S}
 \end{eqnarray}

Let us show that for any $1\le g< i< j\le m$ we have
 \beq{Sgij}
 |\btu(g,i)|+|\btu(j,i)|+|\btu(j,g)|\ge |\btu(g,g)|.
 \eeq
 Start with any $(x,y)\in \btu(g,g)$. 
Let us transform $(x,y)$ into $(\sigma_g(x),\sigma_g(y))$ in three steps,
where we consecutively apply $(\sigma_g,\sigma_i)$,
$(\sigma_j^{-1},\sigma_i^{-1})$, and $(\sigma_j,\sigma_g)$:
 $$
 (x,y)\to(\sigma_g(x),\sigma_i(y))\to (\sigma_j^{-1}(\sigma_g(x)),y)\to
 (\sigma_g(x),\sigma_g(y)). 
 $$
 Since $(x,y,\sigma_g(x),\sigma_g(y))\in \btu$, 
at least one of these three steps changes
adjacency. Depending on the number of the step when this happens, we get
respectively that $(x,y)\in \btu(g,i)$, $(\sigma_j^{-1}(\sigma_g(x)),y)\in \btu(j,i)$,
or $(\sigma_j^{-1}(\sigma_g(x)),y)\in \btu(j,g)$. Conversely, suppose that 
we are given the
resulting conclusion of the form $(u,v)\in\btu(a,b)$ with distinct $a,b\in
\{i,j,g\}$. The pair $(a,b)$ determines the number $k\in \{1,2,3\}$ of the
step. This $k$,
when combined with $(u,v)$, 
easily allows us to reconstruct the ordered pair $(x,y)$. Thus no element in
the left-hand side of~\req{Sgij} is doubly counted. This
proves~\req{Sgij}. 

By~\req{Sgij} (and $|\btu(a,b)|=|\btu(b,a)|$) we conclude that
$S_{\{g,i,j\}}\ge |\btu(g,g)|\ge 2\ell$.
A simple averaging over all choices of $\{i,g,h\}\in {[m]\choose 3}$
implies  that $S_{[m]}\ge 2\ell {m\choose  
2}/{3\choose 2}=\ell m(m-1)/3$. By~\req{delta1S}, we have
 $$
 \delta_1(G,H,P)\ge \frac{2\ell m(m-1)/3 + 2\ell m}{m^2n^2}\ge
 \frac{2\ell}{3n^2} = \frac{\hat\delta_1(G,H)}3,
 $$ 
 finishing the proof of Lemma~\ref{lm:relate}.\epf

\brm The author thanks Alexander Razborov for the remarks that
simplified the original proof of
Lemma~\ref{lm:relate}.\medskip

The interesting problem of finding the best possible constant in
Lemma~\ref{lm:relate} remains open. At the moment, we know only that it is
between $6/5$ (see~\cite[Appendix B]{goldreich+krivelevich+newman+rozenberg:08})
and~$3$.


The situation for the cut-distance is somewhat similar: the discrete version
$\hat\delta_\Box$ of $\delta_\Box$, as defined
by~\cite[Equation~(2.6)]{BCLSV:08}, is not always equal to
the $\delta_\Box$-distance (\cite[Section~5.1]{BCLSV:08}) 
while for any two graphs $G$ and $H$ of the same
order we have
 $$
 \delta_\Box(G,H)\le \hat\delta_\Box(G,H)\le
32(\delta_\Box(G,H))^{1/67}
 $$
 (\cite[Theorem~2.3]{BCLSV:08}). 
 It is open whether $\hat\delta_\Box(G,H)$
can be bounded from above by a linear function of $\delta_\Box(G,H)$,
see e.g~\cite[Page~1830]{BCLSV:08}.

\section{Characterization of Stability}\llabel{stab}

Recall that in the Introduction we defined when an extremal $(f,\C
P)$-problem is stable. Here we give an alternative characterization.
Since stability deals with relating the $\delta_1$ and $\delta_\Box$
distances, it is not surprising that the methods developed by Lov\'asz
and Szegedy~\cite{lovasz+szegedy:08:arxiv} in the context of property
testing apply here.

\bth{stable} Let $\C P$ be an arbitrary graph property with $\C
P_n\not=\emptyset$ for infinitely many $n$ and let $f$ be a
graph parameter. Then the extremal
$(f,\C P)$-problem is stable if and only if $\Lim(f,\C P)$ consists of
a single graphit $[W]$, where moreover
$W\in\WO$ can be chosen to assume values
$0$ and $1$ only.\end{theorem}

The rest of this section is dedicated to proving
Theorem~\ref{th:stable}, in the course of which we observe an interesting
dichotomy result (Theorem~\ref{th:dichotomy}).

We will need the following result, which is a special case 
of~\cite[Lemma~2.2]{lovasz+szegedy:08:arxiv}.

\blm{2.2} Let $W,W_1,W_2,\dots \in \W$ be such that
$\|W_n-W\|_\Box\to 0$ as $n\to\infty$. Let $S\in\C L_{I^2}$.
Then $\int_S W_n\,\dd\lmu \to\int_S
W\,\dd\lmu$ as $n\to\infty$.\end{lemma}

\noindent\textit{Sketch of Proof.} If $S$ is a rectangle, then the
conclusion follows from the definition of the cut-norm. A general
$S\in \C L_{I^2}$ can be approximated within any $\e>0$ by a finite
union of disjoint rectangles, cf Lemma~\ref{lm:approx}.\epf

\bth{dichotomy} Let $W\in\WO$ and let $W_1,W_2,\dots\in\WO$ be an
arbitrary sequence such that $\delta_\Box(W_n,W)\to 0$ as
$n\to\infty$.

If $\lmu(W^{-1}(\{0,1\}))=1$ (that is, $W$ assumes only values $0$
and $1$ a.e.), then the sequence $(W_n)_{n\in\I N}$ is necessarily
convergent to $W$ in the $\delta_1$-distance.
 
If $\lmu(W^{-1}(\{0,1\}))<1$ and each $W_n$ is a.e.\
$\{0,1\}$-valued, then the sequence $(W_n)_{n\in\I N}$ does not
contain any Cauchy subsequence with respect to the
$\delta_1$-distance.
 \end{theorem}
 \bpf Suppose first that $W$ is $\{0,1\}$-valued a.e. Let
$S=W^{-1}(0)\in \C L_{I^2}$. For each $n\in \I N$ choose
$\phi_n\in\Phi_0$ such that $\|W_n^{\phi_n}-W\|_\Box \le
\delta_\Box(W_n,W)+1/n$. Clearly, $\|W_n^{\phi_n}-W\|_\Box$ tends to
0, so by Lemma~\ref{lm:2.2} we have
 \begin{eqnarray*}
 \delta_1(W_n,W)&\le& \|W_n^{\phi_n}-W\|_1 \ =\ \int_SW_n^{\phi_n}\,\dd\lmu +
 \int_{I^2\setminus S} (1-W_n^{\phi_n})\,\dd\lmu \\
 &\to &   \int_S W\,\dd\lmu + \int_{I^2\setminus S} (1-W)\,\dd\lmu\ =\ 0.
 \end{eqnarray*}

Now, suppose that $\lmu(W^{-1}(\{0,1\}))<1$ and that the second part
of the theorem is
false. By choosing a subsequence and relabeling, we can assume that
$(W_n)_{n\in\I N}$ itself is a Cauchy sequence with
$\delta_1(W_m,W_n)\le 1/2^m$ for every $m\le n$. Let $\phi_1:I\to I$ be the
identity map and $U_1=W_1$. Inductively on $n=2,3,\dots$, do the
following. By induction, we assume that we have $U_{n-1}=W_{n-1}^{\phi_{n-1}}$
with $\phi_{n-1}\in \Phi_0$. By Corollary~\ref{cr:inverse},
 $$
 \delta_1(U_{n-1},W_n)=\delta_1(W_{n-1}^{\phi_{n-1}},W_n)=
 \delta_1(W_{n-1},W_n) \le \frac1{2^{n-1}}.
 $$
 Thus there is $\phi_n\in\Phi_0$ such that, letting  $U_n=W^{\phi_n}$,
 we have
 \beq{ind}
 \|U_{n-1}-U_n\|_1\le \frac1{2^{n-2}}.
 \eeq

The sequence $(U_n)_{n\in\I N}$ is Cauchy with respect to
the $\ell_1$-norm: for $m\le n$ we have
 $$
 \| U_n-U_m\|_1\le \sum_{i=m+1}^n \|U_i-U_{i-1}\|_1\le \sum_{i=m+1}^n
   \frac1{2^{i-2}} < \frac{1}{2^{m-2}}. 
 $$
 Since the normed space $L^1$ defined by~\req{L1} is complete
(\cite[Theorem~6.6]{folland:ra}), the sequence $(U_n)_{n\in\I
N}$ has a limit $U\in L^1$:
 \beq{Ulim}
 \lim_{n\to\infty}\| U_n-U\|_1=0.
 \eeq
 We have $\int_{I^2}|U(x,y)-U(y,x)|\,\dd\lmu(x,y)=0$ because it is at
most 
 $$
 2\,\|U-U_n\|_1+ \int_{I^2}|U_n(x,y)-U_n(y,x)|\,\dd\lmu(x,y)=
 2\,\|U-U_n\|_1 \to0.
 $$ Thus $U$ is
symmetric a.e.\ on $I^2$ by e.g.~\cite[Proposition~2.16]{folland:ra}. 
Likewise, $0\le U(x,y)\le 1$ a.e. 
By changing $U$ on a subset of $I^2$ of 
measure zero, we can assume that $U\in\WO$. By the Triangle Inequality,
 $$
 \delta_\Box(U,W)\le \delta_\Box(U,U_n) + \delta_\Box(U_n,W)\le
 \delta_1(U,U_n) + \delta_\Box(W_n^{\phi_n},W).
 $$
 This tends to 0 as $n\to\infty$. Thus $\delta_\Box(U,W)=0$ and by
Theorem~\ref{th:tfae}, $U^\psi=W^\phi$ a.e.\ for some $\psi,\phi\in\Phi$. 
Thus $U$ is not
$\{0,1\}$-valued a.e. 

For $m\in \I N$, let 
 $$
 A_{m}=\{(x,y)\in I^2\mid
1/m< U(x,y)<1-1/m\}.
 $$
 Each $A_m$ is Lebesgue measurable since $U$ is measurable. Also,
$Z=\cup_{m\in \I N} A_m=\{z\in I^2\mid U(z)\not\in\{0,1\}\}$ has positive
measure $c$. By the continuity from below
\cite[Theorem~1.8.c]{folland:ra} of the measure $\lmu$, there is
$m\in \I N$ with
$\lmu(A_{m})>c/2$. Since each $U_n=W_n^{\phi_n}$ is $\{0,1\}$-valued by
assumption, we have $\|U_n-U\|_1\ge c/2m$. This contradicts~\req{Ulim},
and finishes the proof of the lemma.\epf

\brm The first part of Theorem~\ref{th:stable} can also be deduced from~\cite[Lemma~2.9]{lovasz+szegedy:08:arxiv}.\medskip

\begin{corollary}\llabel{cr:dichotomy} Let a
sequence of graphs $G_1,G_2,\dots$ converge in the
$\delta_\Box$-distance to a 
graphit $[W]$. Then the sequence $(G_n)_{n\in\I N}$ converges to $[W]$
in the $\delta_1$-distance if and only if $W$ is $\{0,1\}$-valued
a.e.\qed\end{corollary}

\noindent\textit{Proof of Theorem~\ref{th:stable}:} Suppose first that
the extremal $(f,\C P)$-problem is stable, as defined in
Section~\ref{intro}. Let $[U],[W]\in \Lim(f,\C P)$. Choose witnesses of
this, that is, sequences of almost extremal graphs $(G_{m_i})_{i\in \I
N}$ and $(H_{n_i})_{i\in\I N}$ with $G_{m_i}\to U$ and $H_{n_i}\to W$
in the cut-distance as
$i\to\infty$. By stability, $\delta_1(G_{m_i},H_{n_i})\to 0$. Hence,
 $$
 \delta_\Box(U,W)\le \delta_\Box(U,G_{m_i})+
 \delta_\Box(G_{m_i},H_{n_i}) + \delta_\Box(H_{n_i},W)\le
 \delta_1(G_{m_i},H_{n_i}) +o(1)=o(1).
 $$
 Thus $\delta_\Box(U,W)=0$. Since $[U],[W]\in\Lim(f,\C P)$ were arbitrary, 
the limit set $\Lim(f,\C P)$ consists of a single graphit $[W]$. Since
$(G_{m_i})$ is Cauchy with respect to the $\delta_1$-distance, we
conclude by Theorem~\ref{th:dichotomy} that $W$ is  $\{0,1\}$-valued
a.e., proving one direction of the theorem.

Conversely, suppose that $\Lim(f,\C P)=\{\, [W]\,\}$ for a  $\{0,1\}$-valued
$W\in \WO$. Suppose on the contrary that the extremal problem is not
stable. This implies that there is some $\e>0$ such that for every
$i\in\I N$ there are $m_i,n_i\ge i$, $G_{m_i}\in \C P_{m_i}$,
$H_{n_i}\in \C P_{n_i}$ such that
$f(G_{m_i})\ge \ex_f(m_i,\C P)-1/i$, $f(H_{n_i})\ge \ex_f(n_i,\C P)-1/i$,
and
 \beq{far}
 \delta_1(G_{m_i},H_{n_i})\ge \e.
 \eeq
 By choosing a subsequence and
relabeling, we can additionally assume that for every $i<j$ we have 
$m_i\le n_i<m_j \le n_j$.

By the compactness of $(\CX,\delta_\Box)$ we can find a sequence
$i_1<i_2<\dots$ such that $(G_{m_{i_k}})_{k\in\I N}$ is convergent in
the $\delta_\Box$-distance.
Since $(G_{m_{i_k}})_{k\in\I N}$ is a sequence of almost optimal graphs
with increasing orders, its limit is necessarily $[W]$, the unique
element of $\Lim(f,\C P)$. Likewise, we can find a
subsequence $j_1<j_2<\dots$ of $(i_k)_{k\in\I N}$ such that
the graph sequence $(H_{n_{j_k}})_{j\in\I N}$ 
converges to $[W]$ in $\delta_\Box$. Clearly, the intertwined sequence
$(G_{m_{j_1}},H_{n_{j_1}},G_{m_{j_2}},H_{n_{j_2}},\dots)$ still
converges to $[W]$. By Corollary~\ref{cr:dichotomy}, the last
sequence is Cauchy with respect to the $\delta_1$-distance. This
contradicts~\req{far} and finishes the proof of
Theorem~\ref{th:stable}.\epf

\section{The Erd\H os--Simonovits Stability Theorem}\llabel{es}

In this section, we will prove Theorem~\ref{th:ES}. For this purpose,
we adopt the nice proof of Erd\H os~\cite{erdos:70} that every
$K_{r+1}$-free graph $G$ is \emph{dominated} by some $r$-partite graph
$H$, that is, $V(H)=V(G)$ and $d_H(x)\ge d_G(x)$ for every $x\in
V(G)$, where e.g.\ $d_H(x)$ denotes the degree of $x$ in $H$. In order
to prove this, Erd\H os~\cite{erdos:70} uses induction on $r$ as
follows. The case $r=1$ is trivially true. Let $x$
be a vertex of maximum degree in $G$ and $V'$ be the set of neighbors
of $x$. Then, $G[V']$ is $K_r$-free, so by the induction assumption we can
find an $(r-1)$-partite graph $H'$ that dominates $G[V']$. Let $H$ be
the $r$-partite graph obtained from $H'$ by adding a new part on
$V(G)\setminus V'$. It is not hard to check that $H$ is the required
graph, see~\cite{erdos:70} for details.

Unfortunately, our proof of the graphon version of this
degree-domination result (Theorem~20 in the previous
version~\cite{pikhurko:08:arxiv} of this manuscript) is quite long and
complicated. Later, during a discussion with Peter Keevash, it
was realized that if one is content to prove just Theorem~\ref{th:ES},
then the arguments dealing with graphons can be shortened. Here we present
the shorter proof, referring the interested
reader to~\cite{pikhurko:08:arxiv} for the more general result.

Since we are going to apply the Fubini Theorem a few times, we state
it here. For a function $W: I^2\to\I R$ and $x\in I$, let the
\emph{section functions} $W_x,W^x:I\to\I R$ be defined by
$W_x(y)=W(x,y)$ and $W^x(y)=W(y,x)$. Let $W\Int(x)=\int_I
W_x(y)\,\dd\lmu(y)$ and $W\UInt(x)=\int_I W^x(y)\,\dd\lmu(y)$ (and let
it be arbitrary if the integral is undefined). Clearly, for a
symmetric $W$, we have $W_x=W^x$ and $W\Int=W\UInt$. Since $(I^2,\C
L_{I^2},\lmu_{I^2})$ is not the product $(I,\C L_I,\lmu_I)\times
(I,\C L_I,\lmu_I)$ but its completion, we have to use the Fubini
Theorem for Complete Measures (\cite[Theorem 2.39]{folland:ra}) which
easily follows from the standard Fubini Theorem (\cite[Theorem
2.37.a]{folland:ra}), with the derivation being described in 
\cite[Exercise 2.49]{folland:ra}.

\bth[The Fubini Theorem for the Lebesgue Measure]{fubini}
 If $W\in L^1(I^2,\C L_{I^2},\lmu_{I^2})$, then $W_x,W^x\in
L^1(I,\C L_I,\lmu_I)$ for a.e.\ $x\in I$. Furthermore, $W\Int,W\UInt\in L^1(I,\C
L_I,\lmu_I)$ and
 $$
 \int_{I^2} W(x,y) \,\dd\lmu(x,y) = \int_I W\Int(x)\,\dd\lmu(x) =
 \int_IW\UInt(x)\,\dd\lmu(x).\qed
 $$
 \end{theorem}

Let $W\in \WO$ and $F$ be a graph on $[n]$. We call $W$
\emph{$F$-free} if for every (not necessarily distinct)
$x_1,\dots,x_n\in I$ there is a pair $\{i,j\}\in E(F)$ such that
$W(x_i,x_j)=0$. Equivalently, $W$ is $F$-free if and only if
$W(x,x)=0$ for every $x\in I$ and there is no
homomorphism from $F$ to the infinite (uncountable) graph with vertex
set $I$ in which $x,y$ are connected if $W(x,y)>0$. 

If $W\in\WO$ is $F$-free, then $t(F,W)=0$. The converse is not true:
for example, fix distinct $x_1,\dots,x_n\in I$ and let $W(x,y)=0$ except
$W(x_i,x_j)=1$ for all distinct $i,j\in [n]$. However, please note the
following Lemma~\ref{lm:irl}, which is a rewording of a special case
of a result of Elek and
Szegedy~\cite[Lemma~3.4]{elek+szegedy:08:arxiv}.

\blm[The Infinite Removal Lemma]{irl} For every $W\in \WO$ there is $U\in \WO$ such that $W=U$ a.e.\ and for every graph $F$ either $t(F,U)>0$ or $U$ is
$F$-free.\end{lemma} 

\noindent\textit{Sketch of Proof.} Let $Z$ be the
Lebesgue set of $W$, as defined by~\req{LDT}. Clearly, $Z\subseteq
I^2$ is symmetric. Let $U(x,y)=W(x,y)$ if $(x,y)\in Z$ and $U(x,y)=0$
otherwise. Since $\lmu(Z)=1$, $U=W$ a.e. Also, if $x_1,\dots,x_n$
give an $F$-subgraph in $U$, then there is $c>0$ such that for any
$\{i,j\}\in E(F)$, the measure of
 $$
 \{(x,y)\in (x_i+c,x_i-c)\times
(x_j-c,x_j+c)\mid W(x,y)> W(x_i,x_j)/2>0\}.
 $$
 is, for example, at least $(1-n^{-2})\cdot 4c^2$. It follows that $t(F,W)>0$.\epf

\brm Note that  $W=U$ a.e.\ implies that $t(F,U)=t(F,W)$ for every
graph $F$.\medskip

For the rest of the section, fix an arbitrary family $\C F$
of graphs.  Recall that $\rho(G)=2e(G)/(v(G))^2$
denotes the edge density and $\forb(\C F)$ consists of all $\C F$-free
graphs. For a graphit $[W]$, define $\rho([W])=t(K_2,[W])$. For
convenience, we just write $\rho(W)$. This is compatible with the
previous definition in the sense that for every graph $G$ we have
$\rho(G)=\rho(W_G)$. Define $r$ by~\req{r} and assume that $r\ge 1$. 

Let $\C A$ consist of those graphits $[W]$ that maximize $\rho(W)$
given that $t(F,W)=0$ for every $F\in \C F$.  By the compactness of
$(\C X,\delta_\Box)$ and the continuity of each function $t(F,{-})$,
the maximum is attainable. Denote this maximum value by $a$.

\blm{TuranLimit}
$\Lim(\rho,\forb(\C F))=\C A$.\end{lemma}

\bpf Let $[W]\in \Lim(\rho,\forb(\C F))$. Pick a sequence of almost
extremal graphs $(G_{n_i})$ convergent to $[W]$. Since each $G_{n_i}$
is $\C F$-free, we have $t(F,W)=0$ for each $F\in\C F$ by the first
definition~\req{tFWLimit} of $t(F,W)$. We conclude that $\rho(W)\le a$.

Thus, in order to prove the lemma, it is enough to construct for every
$[W]\in \C A$ a sequence of $\C F$-free graphs $(G_n)_{n\in\I N}$
convergent to $[W]$ with $\rho(G_n)\ge a-o(1)$.  Let $U\in [W]$ be
obtained from $W$ by applying Lemma~\ref{lm:irl}. For each integer $n$
we generate a random graph $G_n$ on $[n]$ as follows. Pick uniformly
at random $n$ elements $x_1,\dots,x_n\in I$ and let a pair $\{i,j\}$
be an edge of $G_n$ with probability $U(x_i,x_j)$, with all
$n+{n\choose 2}$ random choices being mutually independent. With
probability $1$ we have that the sequence $(G_n)$ converges to
$[U]=[W]$ (see \cite[Theorem~4.5]{BCLSV:08} or
\cite[Corollary~2.6]{lovasz+szegedy:06}).  Thus
at least one such sequence $(G_n)$ exists. In particular, we have
$\lim_{n\to\infty} \rho(G_n)= \rho(U)=a$. Also, since $U$ does not
contain any copy of $F\in \C F$, each $G_n$ is (surely) $\C F$-free,
as desired.\epf

\brm The above proof, which is
applicable to many other extremal problems, gives another
justification why it is better not to restrict ourselves to extremal
graphs when defining the limit set $\Lim(f,\C P)$.\medskip

Lemma~\ref{lm:irl} implies that
 \beq{CA}
 \C A =\{\,[W]\mid \forall\, F\in \C F\ F\not\subseteq W \mbox{ and }
 \rho(W) \mbox{ is maximum}\,\}.
 \eeq

Since $W\Int$ is the analytic analog of the degree
sequence, the following lemma can be informally rephrased that
extremal graphons are degree-regular. The combinatorial interpretation of the
proof is that if we have too much discrepancy between degrees in an
almost extremal graph $G$, then by deleting $\e v(G)$ vertices of smaller degree and cloning
$\e v(G)$ vertices of larger degree, we would substantially increase the size of
$G$, which would be a contradiction (provided
we do not create any forbidden subgraph).

\blm{regular} For every $[W]\in \C A$ we have $W\Int(x)=a$ for a.e.\
$x\in I$.\end{lemma}

\bpf The Fubini Theorem implies that if $W=U$ a.e., then $W\Int=U\Int$
a.e. (Indeed, if e.g.\ $W\Int>U\Int$ on a set $X\subseteq I$ 
of positive measure,
then $\int_{X\times I}(W-U)=\int_X (W\Int-U\Int)>0$, a contradiction.)
Hence we can assume by 
Lemma~\ref{lm:irl} that $W$ is $\C F$-free. 

Suppose on the contrary that the lemma is false. Let $X_n=\{x\in I\mid
W\Int(x)\le a-1/n\}$ and $Y_n=\{y\in I\mid W\Int(y)\ge a+1/n\}$. Note
that e.g.\ $\cup_{n\in\I N} X_n=\{x\in I\mid W\Int(x)<a\}$. Since
$\cup_{n\in\I N} (X_n\cup Y_n)$ has positive measure, there is some
$n$ with $\lmu(X_n\cup Y_n)>0$. Assume, for example, that $\lmu(Y_n)$
is positive. By the Fubini Theorem (and $\int_{I^2}W=a$), we conclude
that $\lmu(\cup_{m\in\I N} X_m)>0$. By increasing $n$, assume that
$c=\min(\lmu(X_n),\lmu(Y_n))$ is positive.

Let $\e=\min(c,1/(3n))$. By Lemma~\ref{lm:new}, we can find $\phi\in\Phi_0$ such
that $\lmu(\phi([0,\e])\setminus X_n)=0$ and
$\lmu(\phi([\e,2\e])\setminus Y_n)=0$. Let $U=W^\phi$. Then $U\in [W]$
is still $\C F$-free while $U\Int(x)$ is at most $a-1/n$ (resp.\ at
least $a+1/n$) for a.e.\ $x$ in the interval $[0,\e]$ (resp.\
$[\e,2\e]$). 

For $x\in I$, let $\psi(x)=x$ if $x\ge \e$ and $\psi(x)=x+\e$ if
$x<\e$. Let $V=U^\psi\in\WO$. (Although $\psi$ is not measure preserving,
this definition makes perfect sense.) Note that $V$ is $\C F$-free: if
$x_1,\dots,x_m\in I$ induce a copy of $F$ in $V$, then
$\psi(x_1),\dots,\psi(x_m)$ induce a copy of $F$ in $U$. Moreover,
 $$
 \rho(V)=\int_I V\Int \ge \int_I U\Int-\int_{[0,\e]} U\Int + \int_{[\e,2\e]} U\Int -(2\e)^2
 \ge a - \e(a-1/n)+\e(a+1/n) -4\e^2>a,
 $$ 
  This contradicts the maximality of $a$.\epf

For disjoint measurable sets $A_1,\dots,A_r\subseteq I$, the
\emph{complete $r$-partite graphon} $K_{A_1,\dots,A_r}$ is the simple
function from $I^2$ to $\{0,1\}$ that assumes value $1$ on $\cup_{i\in
[r]}\cup_{j\in [r]\setminus\{i\}} A_i\times A_j$ and $0$ on the
remaining part of $I^2$. (In other words, $W(x,y)=1$ if $x,y$ come
from two different sets $A_i$ and 0 otherwise.) Clearly,
$K_{A_1,\dots,A_r}$ is $K_{r+1}$-free. 

Next, we prove that the graphon problem has the unique solution when
we forbid the clique $K_{r+1}$ only.


\blm{Kr+1} If $\C F=\{K_{r+1}\}$ and $[W]\in\C A$, then there is a
partition $I=A_1\cup\dots\cup A_r$ into sets of measure $1/r$ such
that $W=K_{A_1,\dots,A_r}$ a.e.\end{lemma}

\bpf We use induction on $r$ with the case $r=1$ being trivially true.

Let $r\ge 2$. The $\C F$-free graphon $W_{K_r}$ demonstrates that $a\ge
(r-1)/r$. Let $[W]\in \C A$. Assume that $W$ is $K_{r+1}$-free
by~\req{CA} (that is, by Lemma~\ref{lm:irl}). Pick
$u\in I$ such that $W\Int(u)=a$ which exists  
by Lemma~\ref{lm:regular}. Let $B=\{w\in I\mid W(u,w)>0\}$
and $A_1=I\setminus B$. Let $b=\lmu(B)$. Since $W\le 1$, we have $b\ge
a$. We are free to replace $W$ by $W^\phi$ with any $\phi\in
\Phi_0$; thus  we can assume by  Lemma~\ref{lm:new} that
$\lmu(B\bigtriangleup [0,b])=0$. The graphon $U(x,y)=W(bx,by)$ is
$K_r$-free: if $x_1,\dots,x_r\in I$ induce $K_r$ in $U$, then
$bx_1,\dots,bx_r,u\in I$ induce $K_{r+1}$ in $W$, a contradiction. Note
that, by the Fubini Theorem,
 \beq{a}
 a=\rho(W)= \int_{B^2} W +2\int_{A_1} W\Int - \int_{A_1^2} W =
b^2\rho(U) + 2 (1-b)a -  \int_{A_1^2} W.
 \eeq
 The inductive assumption implies that $\rho(U)\le (r-2)/(r-1)$. Thus
 $$
 \frac{r-1}r\le a\le \frac{r-2}{r-1}\, b^2 + 2(1-b)a\le
 \frac{r-2}{r-1}\, b^2 + 2(1-b)b. 
 $$ 
 Routine algebra implies that $a=b=(r-1)/r$ and all inequalities are in
fact equalities. Thus $W(x,y)=0$ for a.e.\ $(x,y)\in A_1^2$. Since
$W\Int(x)=a=1-\mu(A_1)$ for almost every $x\in A_1$, we have by the Fubini
Theorem that
$W(x,y)=1$ for a.e.\ $(x,y)\in A_1\times B$. Furthermore, by
the uniqueness part of the induction assumption,
$U=K_{B_2,\dots,B_r}$ a.e.\ for some   equitable partiton
$I=B_2\cup\dots\cup B_r$. Letting $A_i=\{bx\mid x\in B_i\}$, we get
$W=K_{A_1,\dots,A_r}$ a.e., as required.\epf

\noindent\emph{Proof of Theorem~\ref{th:ES}:}  By
Theorem~\ref{th:stable}, it suffices to
show that any $[W]\in \Lim(\rho,\forb(\C F))$ we have
$\delta_\Box(W,K_r)=0$.
By Lemma~\ref{lm:TuranLimit} and~\req{CA}, we can assume
that $W$ is $\C F$-free. 

Let us show that $W$ is $K_{r+1}$-free. Suppose on the contrary that
$x_1,\dots,x_{r+1}\in I$ induce $K_{r+1}$ in $W$. Select $F\in \C F$
of chromatic number $\chi(F)=r+1$ and fix a proper coloring $c:V(F)\to
[r+1]$. Then the map $f:V(F)\to I$ with $f(u)=x_{c(u)}$ shows that
$F\subseteq W$, a contradiction.

Tur\'an graphs $T_r(n)$ (or the graphon $W_{K_r}$ and
Lemma~\ref{lm:TuranLimit}) show that $\rho(W)\ge (r-1)/r$. By
Lemma~\ref{lm:Kr+1} we have that $\rho(W)\le\rho(K_r)\le (r-1)/r$. Thus $[W]$ is
extremal for the $(\rho,\forb(\{K_{r+1}\}))$-problem and (again by
Lemma~\ref{lm:Kr+1}) is equal to
$[W_{K_r}]$, as required.\epf

\section*{Acknowledgments}

The author thanks Boris Bukh, Tomasz Gdala, Peter Keevash, Po-Shen
Loh, L\'aszl\'o Lov\'asz, Alexander Razborov, Asaf Shapira, and Benny
Sudakov for helpful discussions and/or comments. Also, the author is
grateful to the anonymous referee for careful reading.\medskip

\renewcommand{\baselinestretch}{1.1}

\end{document}